\documentclass[12pt,a4]{article}
\usepackage{epsfig}
\usepackage{amssymb}
\usepackage{graphicx}
\usepackage{amsmath}
\usepackage{amsfonts}
\usepackage{amsthm}
\usepackage{amscd}
\usepackage{arydshln}
\usepackage{mathtools}
\usepackage{units}
\usepackage{hyperref}
\usepackage[all,cmtip]{xy}
\usepackage{diagxy}
\usepackage{xymatrix}
\usepackage[all]{xy}
\usepackage{fancyhdr}
\usepackage[nottoc,numbib]{tocbibind}
\usepackage{graphicx}
\usepackage{mathrsfs}
\usepackage[perpage]{footmisc}
\usepackage{titlesec}
\usepackage{color}
\usepackage{geometry}
\usepackage{multicol}
\usepackage{slashed}
\usepackage{appendix}
\usepackage{changepage}
\usepackage{bbm}
\usepackage{musixtex}

\usepackage{scrextend}
\deffootnote[1em]{1em}{2em}{\textsuperscript{\thefootnotemark}}

\newtheoremstyle{basic}{11pt}{11pt}{}{}{\bfseries}{.}{0.5em}{}
\newtheoremstyle{proof}{11pt}{11pt}{}{}{\scshape}{:}{0.5em}{}

\newtheorem{prop}{Proposition}

\newtheorem{thm}{Theorem}

\theoremstyle{basic}
\newtheorem{defn}{Definition}

\theoremstyle{basic}
\newtheorem{exa}{Example}

\theoremstyle{basic}

\theoremstyle{proof}
\newtheorem*{prf}{Proof}

\newcommand{\harpoon}{\overset{\rightharpoonup}}

\long\def\comments#1{}

\DeclareMathOperator{\Pull}{Pull}

\makeatletter
\DeclareRobustCommand\bigop[1]{%
  \mathop{\vphantom{\sum}\mathpalette\bigop@{#1}}\slimits@
}
\newcommand{\bigop@}[2]{%
  \vcenter{%
    \sbox\z@{$#1\sum$}%
    \hbox{\resizebox{\ifx#1\displaystyle.9\fi\dimexpr\ht\z@+\dp\z@}{!}{$\m@th#2$}}%
  }%
}
\makeatother

\begin{document}

\pagestyle{plain}

\title{Weak Algebra Bundles and Associator Varieties}

\author{Clarisson Rizzie Canlubo\footnote{\url{crpcanlubo@math.upd.edu.ph}}}

\maketitle

\begin{abstract}
Algebra bundles, in the strict sense, appear in many areas of geometry and physics. However, the structure of an algebra is flexible enough to vary non-trivially over a connected base, giving rise to a structure of a weak algebra bundle. We will show that the notion of a weak algebra bundle is more natural than that of a strict algebra bundle. We will give necessary and sufficient conditions for weak algebra bundles to be locally trivial. The collection of non-trivial associative algebras of a fixed dimension forms a projective variety, called associator varieties. We will show that these varieties play the role the Grassmannians play for principal $O(n)$-bundles.

\


\

\end{abstract}

\

\section{Introduction} \label{S1}

Weak algebra bundles are generalizations of (strict) algebra bundles. They are monoid objects in the category of vector bundles. Algebra bundles appear more frequently in the literature. The exterior bundle and the Clifford bundle are examples. Algebra bundles are examples of weak algebra bundles. In Section \ref{AV}, we look at the varieties of associative algebras of a fixed dimension, the so-called associator varieties. In Section \ref{CWAB}, we will show that weak algebra bundles are more natural than algebra bundles by constructing the so-called classifying weak algebra bundle. In Section \ref{LTWAB}, we will give necessary and sufficient conditions for a weak algebra bundle to be locally trivial, and hence strictness. We will introduce the notion of a differential connection. Existence of a differential connection together with a technical condition guarantee local triviality.

A lot has been written for (associative) algebra bundles. See for example Chidambara-Kiranagi \cite{chi} and Kiranagi-Rajendra \cite{kir}. An almost equal amount of literature has been devoted to Lie algebra bundles. See for example Douady-Lazard \cite{dou}, and the series of papers by Kiranagi et. al. \cite{k2}, \cite{k3}, \cite{k4}, \cite{k5}, \cite{k6}, \cite{k7}, \cite{k8}, and \cite{k9}. In this article, we restrict to finite-rank weak algebra bundles. For the infinite dimensional case, many work has been done. See for example Dadarlat \cite{dad}.

\section{Weak Algebra Bundles}\label{WA}

An \textit{algebra bundle} is a vector bundle in which the fibers are algebras rather than just vector spaces and such that the trivialization maps are algebra isomorphisms. It follows immediately the if the base space is connected then the fiber algebras are mutually isomorphic. For a weak algebra bundle, we do not require that the local trivialization maps are algebra isomorphisms. A \textit{weak algebra bundle} over a space $X$ is a monoid object in the category of real vector bundles $Vec(X)$ over $X$. More precisely, we have the following definition.

\begin{defn}
Let $X$ be a topological space. A \textit{weak algebra bundle} $A\stackrel{p}{\longrightarrow}X$ is a vector bundle together with bundle map $A\otimes A \stackrel{\mu}{\longrightarrow} A$ making the following diagram of bundle maps

\[ \xymatrix{ A\otimes A\otimes A \ar[rr]^-{id \ \otimes \mu} \ar[dd]_-{\mu\otimes \ id} && A\otimes A \ar[dd]^-{\mu}\\
&& \\
A\otimes A \ar[rr]_-{\mu} && A} \]

\noindent commute. If, in addition, there is a bundle map $\mathbbm{1}_{X}\stackrel{\eta}{\longrightarrow}A$, where $\mathbbm{1}_{X}$ is the trivial line bundle over $X$, making the following diagram

\[ \xymatrix{\mathbbm{1}_{X}\otimes A \ar@{=}[rr] \ar[rrdd]_-{\eta\otimes \ id} && A  \ar@{=}[rr] && A\otimes \mathbbm{1}_{X} \ar[lldd]^-{id \ \otimes \eta} \\
&& \\
&& A\otimes A \ar[uu]^-{\mu} &&} \]

\noindent  commute, then $A\stackrel{p}{\longrightarrow}X$ is called \textit{unital}. The map $\mu$ is called the \textit{bundle multiplication} or simply the multiplication while the map $\eta$ is called the \textit{bundle unit} or simply the unit.
\end{defn}

Note that $\Gamma(X,A\otimes A)\cong \Gamma(X,A)\otimes_{C(X)}\Gamma(X,A)$ and $\Gamma(X,\mathbbm{1}_{X})\cong C(X)$ as $C(X)$-bimodules. The global section functor $\Gamma$ induces a multiplication $\mu_{\ast}$ and a unit map $\eta_{\ast}$ on $\Gamma(X,A)$ given by

\[ (\sigma\tau)(x)=\sigma(x)\tau(x) \hspace{.25in} \eta_{\ast}(\alpha)(x)=\eta(\alpha(x))  \]

\noindent for any $\sigma, \tau\in \Gamma(X,A)$, $\alpha\in C(X)$, and $x\in X$. These maps turn $\Gamma(X,A)$ into a unital $C(X)$-algebra. Conversely, a unital $C(X)$-algebra structure on $\Gamma(X,A)$ turns $A\stackrel{p}{\longrightarrow} X$ into a weak algebra bundle. From this equivalence, we see immediately that (strict) algebra bundles are weak algebra bundles.

\section{Associator Varieties}\label{AV}

Despite the name, weak algebra bundles are more natural than (strict) algebra bundles. In this section, we will construct a certain universal weak algebra bundle which describes \textit{all} weak algebra bundles (strict included) of a particular rank. First, let us consider a finite-dimensional vector space $A$ with a chosen basis $\left\{x_{1}, x_{2}, \dots , x_{n} \right\}$. An associative algebra structure on $A$ is completely determined by the structure constants $\alpha_{ij}^{k}$, $1\leqslant i,j,k \leqslant n$ satisfying

\begin{equation}\label{ass}
\sum\limits_{l}\left( \alpha_{ij}^{l}\alpha_{lk}^{m}-\alpha_{il}^{m}\alpha_{jk}^{l}\right) =0.
\end{equation}

\noindent for all $1\leqslant i,j,m,k \leqslant n$. Let $\chi_{n}$ be the variety defined by the $n^4$ equations (\ref{ass}), called the \textit{rank} $n$ \textit{associator variety}. Let $\xi_{n}$ be the quotient of $\chi_{n}$ by the equivalence relation $x\sim y$ if $A_{x}\cong A_{y}$ where $A_{x}$ means the algebra structure on $A$ corresponding to the set of structure constants $x\in \chi_{n}$, equipped with the quotient topology. In general, $\xi_{n}$ is not a variety, just an orbivariety. For $n\in\mathbb{N}$, the \textit{classifying weak algebra bundle} $\widehat{\mathcal{A}}$ of rank $n$ is the weak algebra bundle $\mathcal{A}_{n}\stackrel{p}{\longrightarrow}\xi_{n}$ such that $p^{-1}(x)=A_{x}$.

Let us look more closely to the associator varieties $\chi_{n}$. For the purpose of this section, we let $\chi_{n}$ be the variety given by the system of equations \ref{ass} with the all-zero solution removed. This makes $\chi_{n}$ a projective variety. Before we go into the analysis of these varieties, let us give some explicit points.

\begin{exa}
For any $1\leqslant i,j,k \leqslant n$, let
\[ \alpha_{ij}^{k}=\left\{ \begin{array}{lll}
  1   & & k=i+j \ \text{mod} \ n \\
  0   & & \text{otherwise}
\end{array}\right.\]

\noindent The algebra $A_{x}$, where $x=(\alpha_{ij}^{k})$, is the truncated polynomial algebra.
\end{exa}

\begin{exa}
For any $1 \leqslant i,j,k \leqslant n$, let

\[ \alpha_{ij}^{k}=g(k)h(i)h(j) \]

\noindent where $g,h:\left\{1,2,\dots,n \right\}\longrightarrow \mathbb{C}$ are arbitrary functions. A particular example is when $g(k)=k$ and $h(j)=e^{\pi i j}$.
\end{exa}

The associator varieties fit into a natural sequence. The variety $\chi_{n-1}$ is the intersection of $\chi_{n-1}$ with the varieties $\alpha_{in}^{k}=\alpha_{ni}^{k}=0$, $1\leqslant i,k\leqslant n$. This gives an inclusion of varieties $\chi_{n-1}\stackrel{i_{n-1}}{\longrightarrow}\chi_{n}$. This inclusion is natural in the sense of the following proposition.

\begin{prop}\label{P1}
The tautological weak algebra bundle $\mathcal{W}_{n-1}\stackrel{r_{n-1}}{\longrightarrow}\chi_{n-1}$ is the pullback of the tautological weak algebra bundle $\mathcal{W}_{n}\stackrel{r_{n}}{\longrightarrow}\chi_{n}$ along the map $i_{n-1}$ defined above.
\end{prop}

\begin{prf}
The proposition follows directly from the following pullback cube

\[ \xymatrix{ \mathcal{W}_{n-1} \ar[rr] \ar[rd] \ar[dd] && \mathcal{W}_{n} \ar'[d][dd] \ar[rd] & \\
& \mathcal{A}_{n-1} \ar[dd] \ar[rr] && \mathcal{A}_{n} \ar[dd] \\
\chi_{n-1} \ar'[r][rr] \ar[rd]  && \chi_{n} \ar[rd] & \\
& \xi_{n-1} \ar@{-->}[ru]_-{j} \ar[rr]_-{\hat{\imath}} && \xi_{n} } \]

\noindent where the map $j$ is the map induced by the universality of the quotient $\xi_{n-1}$ and the map $\hat{\imath}$ is the composition of the $j$ and the quotient map $\chi_{n}\longrightarrow \xi_{n}$. $\blacksquare$
\end{prf}

By Proposition \ref{P1}, we have several stratifications, given by the horizontal maps in the complex of spaces below.

\[ \xymatrix{ \mathcal{W}_{n-1} \ar[rr] \ar[rd] \ar[dd] && \mathcal{W}_{n} \ar[rr] \ar[rd] \ar'[d][dd] &&  \mathcal{W}_{n+1} \ar[rr] \ar[rd] \ar'[d][dd] && \mathcal{W}_{n+2} \ar'[d][dd] \ar[rd] & \\
& \mathcal{A}_{n-1} \ar[rr] \ar[dd] && \mathcal{A}_{n} \ar[rr] \ar[dd] && \mathcal{A}_{n+1} \ar[dd] \ar[rr] && \mathcal{A}_{n+2} \ar[dd] \\
\chi_{n-1} \ar'[r][rr] \ar[rd] && \chi_{n} \ar'[r][rr] \ar[rd] && \chi_{n+1} \ar'[r][rr] \ar[rd]  && \chi_{n+2} \ar[rd] & \\
& \xi_{n-1} \ar[rr] && \xi_{n} \ar[rr] && \xi_{n+1} \ar[rr] && \xi_{n+2} }\]

Using maps above, we can define the spaces $\chi_{\infty}=\lim\limits_{\longrightarrow}\chi_{n}$. Similarly, we define  $\xi_{\infty}$, $\mathcal{W}_{\infty}$ and $\mathcal{A}_{\infty}$ as direct limits of the obvious sequence of spaces and maps. Then, it is immediate to check that there are weak algebra bundles $\mathcal{W}_{\infty}\stackrel{r}{\longrightarrow}\chi_{\infty}$ and $\mathcal{A}_{\infty}\stackrel{p}{\longrightarrow}\xi_{\infty}$. Moreover, these fit into a pullback square

\[ \xymatrix{ \mathcal{W}_{\infty} \ar[rr] \ar[dd]_-{r} && \mathcal{A}_{\infty} \ar[dd]^-{p} \\
&& \\
\chi_{\infty} \ar[rr] && \xi_{\infty}} \]

\noindent where $\chi_{\infty}\longrightarrow \xi_{\infty}$ is the direct limit of the quotient maps $\chi_{n}\longrightarrow \xi_{n}$.

Let us end this section by looking at the tangent spaces of associator varieties. The following computation can be found in \cite{shafarevich1}. Consider a point $\alpha\in\chi_{n}$ given by $\left\{\alpha_{ij}^{k} | 1\leqslant i,j,k\leqslant n\right\}$ in $\chi_{n}$. Then tangent vectors $v=\left\{v_{ij}^{k} | 1\leqslant i,j,k\leqslant n\right\}$ to $\chi_{n}$ at the point $\alpha$ satisfy the equation

\begin{equation}\label{tan}
\sum_{l} \left( \alpha_{ij}^{l}v_{lk}^{m}+\alpha_{lk}^{m}v_{ij}^{l}-\alpha_{il}^{m}v_{jk}^{l}-\alpha_{jk}^{l}v_{il}^{m} \right) =0
\end{equation}

\noindent Let $x_{1},\dots,x_{n}$ be a basis for $A$ giving the structure constants $\left\{\alpha_{ij}^{k} | 1\leqslant i,j,k\leqslant n\right\}$. The bilinear function $f_{v}:A\times A\longrightarrow A$ given by $f(x_{i},x_{j})=\sum_{k}v_{ij}^{k}x_{k}$ satisfies the relation

\begin{equation}\label{cocycle}
    xf(y,z) - f(xy,z) + f(x,yz) -f(x,y)z = 0
\end{equation}
\

\noindent for all $x,y,z\in A$. Condition \ref{cocycle} is called the cocycle condition. In particular, bilinear functions $f$ satisfying \ref{cocycle} are called $2$-cocycles. The set of all $2$-cocycles in $A$, denoted by $Z^{2}(A)$ is a vector space. The map $v\mapsto f_{v}$ defines a linear isomorphism between the tangent space to $\chi_{n}$ at the point $\alpha$ and the space of $2$-cocycles.

\section{Classifying Weak Algebra Bundles}\label{CWAB}

In this section, we will show that the orbivariety $\xi_{n}$ is the analogue the Grassmannians play for principal $O(n)$-bundles. Let us first recall the the Grassmannian variety $G(n,\mathbb{R}^{\infty})$ is the \textit{classifying space} for $O(n)$. This in particular means that for any principal $O(n)$-bundle $Y\stackrel{p}{\longrightarrow}X$ there is a continuous map $X\stackrel{f}{\longrightarrow} G(n,\mathbb{R}^{\infty})$ such that $Y$ is the pull-back of $EO(n)$ along $f$. Here, $EO(n)$ is the total space of the \textit{universal bundle} over $BG=G(n,\mathbb{R}^{\infty})$. See \cite{may} for more details. The following theorem jusitifies the name of the bundle $\mathcal{A}_{n}\stackrel{p}{\longrightarrow}\xi_{n}$.

\begin{thm}\label{class}
Let $\mathcal{A}_{n}\stackrel{p}{\longrightarrow}\xi_{n}$ be the classifying weak algebra bundle of rank $n$. Let $\mathcal{B}\stackrel{q}{\longrightarrow}X$ be a weak algebra bundle of rank $n$. Then there is a continuous map $X\stackrel{f}{\longrightarrow}\xi_{n}$ such that $\mathcal{B}\cong f^{\ast}\mathcal{A}_{n}$ as weak algebra bundles, i.e. a map $f$ which makes the following diagram

\[ \xymatrix{ \mathcal{B} \ar[rr] \ar[dd]_-{q} && \mathcal{A}_{n} \ar[dd]^-{p} \\
&& \\
X \ar[rr]_-{f} && \xi_{n}} \]

\noindent a pullback diagram.

\end{thm}

\begin{prf}
For a weak algebra bundle $\mathcal{B}\stackrel{q}{\longrightarrow}X$, define the map $X\stackrel{f}{\longrightarrow}\xi_{n}$ that sends $x\in X$ to $y\in \xi_{n}$ if $q^{-1}(x)\cong p^{-1}(y)$. Let $\left\{ U_{\alpha}|\alpha\in I\right\}$ be an open cover of $X$ trivializing $\mathcal{B}\stackrel{q}{\longrightarrow}X$. Then $q^{-1}(U_{\alpha})\cong U_{\alpha}\times B$ where $B$ is the underlying vector space of the typical fiber of $q$. Let $x_{1},\dots,x_{n}$ be a basis for $B$ and $\gamma_{ij}^{k}(x)$ be the structure functions of $B_{x}$ for $x\in U_{\alpha}$. Then, the continuous functions $\gamma_{ij}^{k}:U_{\alpha}\longrightarrow \mathbb{R}$ satisfy the equations \ref{ass}. This defines a continuous map $U_{\alpha}\stackrel{f}{\longrightarrow}\xi_{n}$ sending $x\in U_{\alpha}$ to $(\gamma_{ij}^{k}(x)|1\leqslant i,j,k \leqslant n)\in \xi_{n}$. By Proposition \ref{para}, these functions extend globally to a continuous function $X\stackrel{f}{\longrightarrow}\xi_{n}$. Finally, it is straightforward to check that, indeed, $\mathcal{B}\cong f^{\ast}\mathcal{A}_{n}$. $\blacksquare$
\end{prf}

If $\mathcal{B}\stackrel{q}{\longrightarrow}X$ is a (strict) algebra bundle the map $X\stackrel{f}{\longrightarrow}\xi_{n}$ asserted by Theorem \ref{class} is just the constant map, sending every point of $X$ to the unique point $y\in\xi_{n}$ such that the typical fiber of $\mathcal{B}\stackrel{q}{\longrightarrow}X$ is isomorphic to $p^{-1}(y)$.

In view of Proposition \ref{para}, the non-triviality of the underlying vector bundle of a rank $n$ weak algebra bundle $\mathcal{B}\stackrel{q}{\longrightarrow}X$ is determined by the homotopy type of the function $f$ asserted by Theorem \ref{class} but not completely so. For example, a strict algebra bundle $\mathcal{B}\stackrel{q}{\longrightarrow}X$ may have a non-trivial underlying vector bundle and yet the associated function $f$ is constant and hence, homotopically trivial.

Let $\mathcal{W}_{n}\stackrel{r}{\longrightarrow}\chi_{n}$ be the pullback of the rank $n$ classifying weak algebra bundle along the quotient map $\chi_{n}\stackrel{f}{\longrightarrow}\xi_{n}$. We call this bundle the \textit{tautological weak algebra bundle}. The following theorem illustrates that for most purposes we can use the tautological weak algebra bundle in place of the classifying weak algebra bundle. Denote by $\Pull{\left(\mathcal{B},X\right)}$ the set of all weak algebra bundles that are pullbacks of $\mathcal{B}\stackrel{f}{\longrightarrow}X$.

\begin{thm}
There is a canonical bijection between $\Pull{\left(\mathcal{A}_{n},\xi_{n}\right)}$ and $\Pull{\left(\mathcal{W}_{n},\chi_{n}\right)}$.
\end{thm}

\begin{prf}
If $\mathscr{B} \longrightarrow Y$ is a pullback of the classifying weak algebra bundle, then pulling-back maps along appropriate maps as illustrated by the left cube below gives a weak algebra bundle $\mathscr{C}\longrightarrow q^{\ast}Y$ where $q$ is the quotient map.

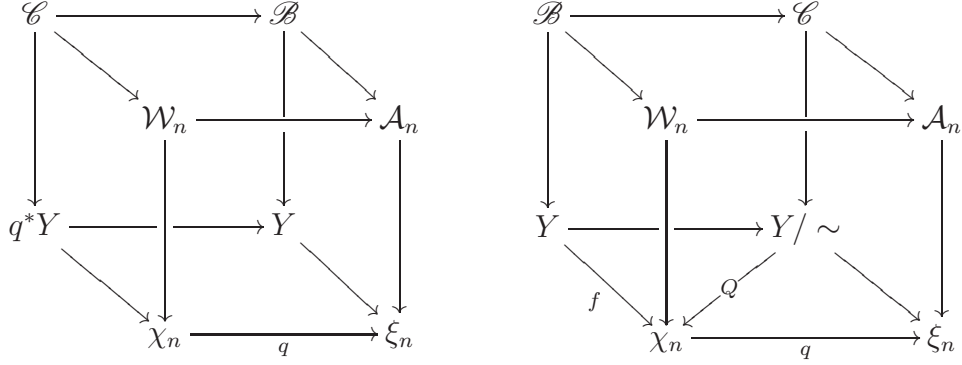
\[ \xymatrix{ \mathscr{C} \ar[rr] \ar[rd] \ar[dd] && \mathscr{B} \ar'[d][dd] \ar[rd] & \\
& \mathcal{W}_{n} \ar[dd] \ar[rr] && \mathcal{A}_{n} \ar[dd] \\
q^{\ast}Y \ar'[r][rr] \ar[rd]  && Y \ar[rd] & \\
& \chi_{n} \ar[rr]_-{q} && \xi_{n} } \hspace{.5in} \xymatrix{ \mathscr{B} \ar[rr] \ar[rd] \ar[dd] && \mathscr{C} \ar'[d][dd] \ar[rd] & \\
& \mathcal{W}_{n} \ar[dd] \ar[rr] && \mathcal{A}_{n} \ar[dd] \\
Y \ar'[r][rr] \ar[rd]_-{f}  && Y/\sim \ar[ld]|-{Q} \ar[rd] & \\
& \chi_{n} \ar[rr]_-{q} && \xi_{n}    } \]

Conversely, let $\mathscr{B} \longrightarrow Y$ be a pullback of the tautological weak algebra bundle along a map $Y\stackrel{f}{\longrightarrow}\chi_{n}$. Let $\sim$ be the equivalence relation on $Y$ defined as $y\sim y'$ if $f(y)=f(y')$. Let $Q$ be the map induced by the universality of the quotient $Y/\sim$ and let $Y/\sim \longrightarrow \xi_{n}$ be the composition of $Q$ and the quotient map $q$. Pulling back maps along appropriate maps according to the right cube above gives a weak algebra bundle $\mathscr{C}\longrightarrow Y/\sim$ that is a pullback of the classifying weak algebra bundle. $\blacksquare$
\end{prf}

The advantage of working with $\mathcal{W}_{n}\stackrel{r}{\longrightarrow}\chi_{n}$ is the fact that $\chi_{n}$ is a variety and $\mathcal{W}_{n}\stackrel{r}{\longrightarrow}\chi_{n}$ is a regular vector bundle. Let us end this section by a triviality statement regarding the classifying weak algebra bundles.

\begin{prop}\label{para}
As vector bundles, $\mathcal{A}_{n}\stackrel{p}{\longrightarrow}\xi_{n}$ are parallelizable for all $n\in\mathbb{N}$.
\end{prop}

\noindent For $i=1,\dots,n$, the sections $\xi_{n}\stackrel{\sigma_{i}}{\longrightarrow}\mathcal{A}_{n}$, $x\mapsto x_{i}$ gives a set of pointwise linearly independent set of $n$ sections. This illustrates parallelizability.

\section{Local Triviality of Weak Algebra Bundles}\label{LTWAB}

In this section, we give necessary and sufficient conditions for a weak algebra bundle to be a strict algebra bundle. For this purpose, we will specialize in the smooth case. Let $X$ be a connected smooth manifold.

\begin{defn}\label{D4.4.1}
Let $E\twoheadrightarrow X$ be a smooth vector bundle such that the fibers are algebras whose multiplications depend on $x\in X$ continuously. A \textit{differential connection} $\nabla$ on $E$ is a smooth connection such that for any vector field $\nu$ on $X$, we have

\[ \nabla_{\nu}(\sigma_{1}\sigma_{2})=\sigma_{1}\nabla_{\nu}(\sigma_{2})+\nabla_{\nu}(\sigma_{1})\sigma_{2} \]

\noindent for any sections $\sigma_{1},\sigma_{2}\in \Gamma(X,E)$. $\Box$
\end{defn}

Surprisingly, existence of such connections is a sufficient condition for the fiber algebras to be isomorphic. For a necessary condition, one needs a stronger assumption than just having isomorphic fiber algebras. We will formalize these statements in the next two propositions.

\begin{prop}\label{4.4.4}
If $E$ has a differential connection $\nabla$ then the fiber algebras of $E\twoheadrightarrow X$ are all isomorphic.
\end{prop}

\begin{prf}
Assume $E$ has a differential connection $\nabla$. Let $x,y\in X$ and let $\gamma:I\longrightarrow X$ be a (piecewise) smooth path in $X$ with $\gamma(0)=x$ and $\gamma(1)=y$. Using the connection $\nabla$, we have a parallel transport map

\[ \Phi(\gamma)^{y}_{x}:E_{x}\longrightarrow E_{y} \]

\noindent which is a linear isomorphism. Thus, all we have to show is that $\Phi(\gamma)^{y}_{x}$ is multiplicative. Given $b_{1},b_{2}\in E_{x}$, there are unique smooth sections $\sigma_{1}$ and $\sigma_{2}$ of $E$ along $\gamma$ such that $\nabla_{\harpoon \gamma}\sigma_{1}=\nabla_{\harpoon \gamma}\sigma_{2}=0$ and $\sigma_{1}(x)=b_{1}$ and $\sigma_{2}(x)=b_{2}$. Here, $\harpoon \gamma$ denotes the smooth tangent vector field of $\gamma$. Note that the product $\sigma_{1}\sigma_{2}$ is the unique smooth section of $\xymatrix{E \ar@{->>}[r] & X}$ along $\gamma$ such that $\left(\sigma_{1}\sigma_{2}\right)(x)=\sigma_{1}(x)\sigma_{2}(x)=b_{1}b_{2}$ and

\[ \nabla_{\harpoon \gamma}(\sigma_{1}\sigma_{2})= \sigma_{1}\nabla_{\harpoon \gamma}(\sigma_{2})+\nabla_{\harpoon \gamma}(\sigma_{1})\sigma_{2} =0. \]

\noindent Thus, by definition of the parallel transport map $\Phi\left(\gamma\right)_{x}^{y}$ we have

\[ \Phi\left(\gamma\right)_{x}^{y}(b_{1}b_{2})=\left(\sigma_{1}\sigma_{2}\right)(y)=\sigma_{1}(y)\sigma_{2}(y)=\Phi\left(\gamma\right)_{x}^{y}(b_{1})\Phi\left(\gamma\right)_{x}^{y}(b_{2}) \]

\noindent which shows that $\Phi\left(\gamma\right)_{x}^{y}$ is multiplicative. $\blacksquare$
\end{prf}

A strong converse of the above proposition, where the isomorphisms among fibers satisfy some coherence conditions, holds. By a coherent collection

\[ \mathscr{P}=\left\{\Phi(\gamma)_{x}^{y}:E_{x}\longrightarrow E_{y}|\forall x,y\in X, \gamma:I\longrightarrow X \ smooth\right\} \]

\noindent of isomorphisms among fibers of $E\twoheadrightarrow X$, we mean a collection satisfying

\begin{enumerate}
\item[(a)] $\Phi(\gamma)^{x}_{x}=id$,
\item[(b)] $\Phi(\gamma)^{y}_{u}\circ\Phi(\gamma)^{u}_{x}=\Phi(\gamma)^{y}_{x}$,
\item[(c)] and $\Phi$ depends smoothly on $\gamma$, $y$ and $x$.
\end{enumerate}

\noindent We then have the following proposition.

\begin{prop}\label{P4.4.5}
A coherent collection $\mathscr{P}$ of algebra isomorphisms on fibers of\\ $E\twoheadrightarrow X$ gives a differential connection $\nabla$ on $E$.
\end{prop}

\begin{prf}

\[ \nabla_{\nu}(\sigma)=\lim\limits_{t\rightarrow 0} \dfrac{\Phi(\gamma)^{x}_{\gamma(t)}\sigma(\gamma(t))-\sigma(x)}{t}=\left.\dfrac{d}{dt}\right|_{t=0}\Phi(\gamma)^{x}_{\gamma(t)}\sigma(\gamma(t)) \]

\noindent for any $\sigma\in B$, $x=\gamma(0)$, and $\nu=\gamma^{'}(0)$. That $\nabla$ is a differential connection follows from the multiplicativity of $\Phi(\gamma)^{y}_{x}$ and the Leibniz property of $\left.\dfrac{d}{dt}\right|_{t=0}$. $\blacksquare$
\end{prf}

\subsubsection*{Acknowledgement}
\footnotesize{I would like to thank Ryszard Nest for supervising me through this project as part of my PhD thesis and for all the stimulating discussions and the Center for Symmetry and Deformation, University of Copenhagen for providing a conducive environment for my studies. I would also like to thank NSRI for funding this project through NSRI-MAT-18-1-02.}

\hspace{1in}

\noindent\textsc{Clarisson Rizzie P. Canlubo}\\
Institute of Mathematics \& \\
National Science Research Institute\\
University of the Philippines$-$Diliman\\
Quezon City, Philippines 1101\\
\url{crpcanlubo@math.upd.edu.ph}

\end{document}